\documentclass[12pt]{article}
\usepackage{amsmath}
\usepackage{amssymb}

\title{Semidomains and Metabelian Product of Metabelian Lie Algebras}
\author{\textsf{Evelina Yu. Daniyarova}\thanks{Supported by RFFI grant
N02-01-00192.} \and \textsf{Ilia V. Kazachkov}\thanks{Supported
by 'Universitety Rossii' grant.} \and \textsf{Vladimir N.
Remeslennikov}$^{\dag}$}

\newtheorem{defn}{Definition}
\newtheorem{thm}{Theorem}
\newtheorem{cor}{Corollary}
\newtheorem{lem}{Lemma}
\newtheorem{prop}{Proposition}

\newtheorem{ex}{Example}

\newcommand{\Fit}{\texttt{Fit}}
\newcommand{\proof}{\paragraph{Proof.}}
\newcommand{\beq}{\begin{equation}}
\newcommand{\eeq}{\end{equation}}

\begin{document}

\maketitle

\section{Introduction}
\label{sec:1}

This paper is the third in a series of papers, the aim of which is to construct algebraic geometry over metabelian Lie
algebras. The foundations of this theory were layed in the papers \cite{AG1, AG2}.

Let $A$ be a metabelian Lie algebra over a field $k$; $S(x)=0$ be
a system of equations over $A$ with the set of variables
$X=\left\{ x_1, \dots, x_n \right\}$ and $V(S)$ be the algebraic
set from $A^n$ defined by system $S$. Analysis of proofs of
theorems from \cite{AG2}, that establish the structure of
algebraic sets $V(S)$ and the structure of their coordinate
algebras $\Gamma (S)$ shows that all major calculations are done
in $A[X]=A\ast F(X)$. Here $F(X)$ stands for the free metabelian
Lie algebra and '$\ast$' denotes the metabelian product. This
argument explains our interest to the problem of investigation of
structure of metabelian product of Lie algebras. The solution to
this problem is given by Theorems \ref{thm:A*B1}, \ref{thm:A*B2}
and \ref{thm:stA*B}.

The notion of a \emph{domain} for groups (the groups that have no
zero divisors) was introduced in \cite{mr}, where the authors
point out its importance for estimation of some certain criterions
of irreducibility of algebraic sets. Every Lie algebra or group
that possesses  non-trivial abelian ideals (subgroups) has zero
divisors, thus in the categories of Lie algebras and metabelian
Lie algebras the notion of a domain fails (in the case of
metabelian Lie algebras the commutant is an abelian ideal). These
circumstances make us to introduce the notions of a
\emph{semidomain} and a \emph{strict semidomain} and investigate
their properties in order to apply obtained results to algebraic
geometry.

\section{Preliminaries}

Below we provide the reader with a brief overview of some
auxiliary facts on Lie algebra (the basics for Lie algebras theory
can be found in \cite{Bah}).

Recall, that a Lie algebra $A$ over a field $k$ is termed
metabelian if and only if the following universal axiom holds:
\begin{itemize}
    \item $(a \circ b) \circ (c \circ d) = 0$.
\end{itemize}

By $a \circ b$ or $ab$ we shall denote the product of elements
from algebra $A$. \emph{Left normed products} $a_1 a_2 a_3 \cdots
a_n$, of $a_1 ,a_2 ,a_3 ,\ldots,a_n \in A$ are defined as
$$
(\ldots((a_1  \circ a_2 ) \circ a_3 ) \circ \ldots) \circ a_n.
$$
We term such products by\emph{ left normed words or monomials} of
the \emph{degree} or of the \emph{length}  $n$.

It is well-known that every monomial of the length $l$  from
letters $a_1, \ldots, a_n$ can be written as a linear combination
of left normed monomials of length $l$ from the same set of
letters.

Let $A$ be an arbitrary Lie algebra over $k$ and let $\left< x
\right>$ denote the principal ideal generated by $x \in A$. In
\cite {Bah} it is shown, that  $\left< x \right>$ is a $k$-linear
span of left normed words from $A$, that begin with $x$:
$$
x, xa_1, xb_1b_2, \ldots, xc_1c_2 \cdots c_n, \ldots
$$

Ideal generated by $\left\{ a \circ b |a, b \in A \right\}$ is
called the \emph{commutant} and is denoted by $A^2$. The
\emph{Fitting's radical} of the Lie algebra $A$ is an ideal,
generated by the set of all elements from nilpotent ideals of $A$.
We'll denote the Fitting's radical of the algebra $A$ by {\rm
$\Fit(A)$}. The Fitting's radical can be characterized as follows:
an element $x$ from $A$ is an element of {\rm $\Fit(A)$} if and
only if the ideal $\left< x \right>$ is nilpotent (see
\cite{AG1}).

Let $F$ be a free algebra in the variety $\mathcal{A}^2$ of all
metabelian Lie algebras.  And let $\left\{a_{\alpha}| \alpha \in
\Lambda \right\}$ be a free base for $F$, here $\Lambda$ is a
totally ordered set. Then left normed monomials $a_{i_1 } a_{i_2 }
\ldots a_{i_m }$ are called \emph{normalised} if and only if the
following condition is satisfied:
$$
i_1  > i_2  \le i_3  \le \ldots \le i_m.
$$
In \cite{Art} it is proven that the set of all normalised words
form a linear basis of $F$.

Let $A$ be an arbitrary metabelian Lie algebra over $k$. Then the
commutant $A^2$ admits a structure of an $R$-module, here $R=k[X]$
is the ring of polynomials over $k$ and $X$ is a linear basis of
vector space $\left. A \right/A^2$ (see \cite{AG1} for details).

\section{Metabelian Products of Metabelian Lie Algebras} \label{sec:prod}

Let $A$ and $B$ be two metabelian Lie $k$-algebras. Throuhgout
this paper the field $k$ is fixed, thus throughout by
$\mathcal{A}^2$ we'll assume the variety of all metabelian Lie
$k$-algebras. Via $A\ast B$ we denote the product of $A$ and $B$
in the variety $\mathcal{A}^2$. Recall the definition of the
algebra $A\ast B$.

\begin{defn}[\cite{Bah}] \label{def:A*B}
Let
$$
A=\left<z_i, i \in I_A|R_A \right>_{\mathcal{A}^2}, \ B=\left<z_i,
i \in I_B|R_B \right>_{\mathcal{A}^2}
$$
be the presentations of respective algebras in the variety
${\mathcal{A}^2}$. Assume that $I_A \cap I_B= \emptyset$. Then the
algebra $A\ast B$ is defined  by the presentation
$$
A\ast B=\left<z_i, i \in I_A \cup I_B| R_A \cup R_B
\right>_{\mathcal{A}^2}.
$$
\end{defn}

\medskip

We also give the categorial definition of $A\ast B$, which is
equivalent to Definition \ref{def:A*B}.

\begin{defn}[Categorial Definition] \label{def:A*B2}
A metabelian Lie algebra $C \in {\mathcal{A}^2}$ is termed
metabelian Lie product of metabelian Lie algebras $A$ and $B$ if
it satisfies two following conditions
\begin{enumerate}
    \item There exist two injections $i_A: A \rightarrow C$ and $i_B: B \rightarrow
    C$ such that the algebra $C$ is generated by the images $i_A(A)$ and
    $i_B(B)$, $C=\left<i_A(A),i_B(B) \right> _{\mathcal{A}^2}$;
    \item For every metabelian Lie algebra $H$ over the field $k$ and
    every pair of injections $\psi_A:A \rightarrow H$ and $\psi_B:B \rightarrow
    H$ there exists a homomorphism $\phi :C \rightarrow H$ such that \/ $\phi i_A=
    \psi_A$ and \/ $\phi i_B= \psi_B$.
\end{enumerate}
\end{defn}

\medskip

In this section we shall investigate the structure of $A\ast B$
and derive that {\rm $\Fit(A\ast B)=(A\ast B)^2$}.

Let $\overline{A}$ and $\overline{B}$ be the vector spaces over
$k$, $\overline{A}={\raise0.7ex\hbox{$A$} \!\mathord{\left/
 {\vphantom {A^2}}\right.\kern-\nulldelimiterspace}
\!\lower0.7ex\hbox{$A^2$}}$ and
$\overline{B}={\raise0.7ex\hbox{$B$} \!\mathord{\left/
 {\vphantom {B^2}}\right.\kern-\nulldelimiterspace}
\!\lower0.7ex\hbox{$B^2$}}$. Assume that $X= \left\{ x_i| i\in
I\right\}$ and $Y= \left\{ y_j| j\in J\right\}$ are the linear
bases of $\overline{A}$ and $\overline{B}$ correspondingly and
that the sets $I$ and $J$ are totally ordered. For the sake of
convenience, by the elements of the sets $X$ and $Y$ we also
denote their fixed preimages in $A$ and $B$.

As we have already mentioned above  the commutant $A^2$ admits the
structure of a module over $R_X=k\left[X\right]$ and the commutant
$B^2$ the structure of
 a module over $R_Y=k\left[Y\right]$. The
$k$-vector space ${\raise0.7ex\hbox{$A\ast B$} \!\mathord{\left/
 {\vphantom {(A\ast B)^2}}\right.\kern-\nulldelimiterspace}
\!\lower0.7ex\hbox{$(A\ast B)^2$}}$ is generated by the set $X
\cup Y$. This implies that $(A\ast B)^2$ can be treated as a
module over the ring of polynomials $R_{X \cup Y}= k\left[X \cup
Y\right]$.

We estimate the structure of the algebra $A\ast B$ in two steps.
\begin{description}
  \item[Step 1] Define a 3-tuple of $R_{X \cup Y}$-modules $M_0$,
    $M_1$ and $M_2$. With the help of these modules we construct a new
    $R_{X \cup Y}$-module $M$.
  \item[Step 2] By the module $M$ we construct a
    metabelian Lie algebra $C$ and prove that $C$ is isomorphic
    to $A\ast B$. This isomorphism yields to quite satisfactory
    structural results for $A\ast B$.
\end{description}

Extend the modules $A^2$ and $B^2$ ($R_X$- and $R_Y$- modules
respectively) to $R_{X \cup Y}$-modules by setting:
$$
M_1=A^2 \otimes_{R_X} R_{X \cup Y}, \ M_2=B^2 \otimes_{R_Y} R_{X
\cup Y},
$$
where here `$\otimes_{R_X}$' is the sign of tensor product over
$R_X$ and `$\otimes_{R_Y}$' stands for  tensor product over $R_Y$.
Let $M_0$ be a module isomorphic to the ideal of $R_{X \cup Y}$
generated by $\left\{ x_iy_j | i \in I, j \in J \right\}$. In
particular, $M_0$ is torsion free. Let
$$
A^2= \left< Z_1| R_1 \right>; \quad B^2=\left< Z_2| R_2 \right>, \
Z_1 \cap Z_2 = \emptyset
$$
be the presentations of $R_X$-module $A^2$ and $R_Y$-module $B^2$.
To define the module $M$ by generators and relations we denote
some of the generators by a pair of letters. Let $XY=\left\{
w_{ij} \equiv x_iy_j | i \in I, j \in J \right\}$, $Z = Z_1 \cup
Z_2$.
 We set
\begin{equation} \label{eq:M}
  M=\left<XY \cup Z | R_1 \cup R_2 \cup S \right>,
\end{equation}
where here
$$
S=\left\{
\begin{array}{ll}
    (x_iy_{j_1})\cdot y_{j_2}= (x_iy_{j_2})\cdot y_{j_1} +( y_{j_2}y_{j_1})\cdot x_i  ; j_1>j_2; \\
    (x_{i_1}y_j)\cdot x_{i_2}= ( x_{i_2}y_j)\cdot x_{i_1} + (x_{i_1}x_{i_2})\cdot y_j ; i_1>i_2. \\
\end{array}
\right.
$$
The symbol `$\cdot$'  denotes the module action  and
$y_{j_2}y_{j_1}$, $x_{i_1}x_{i_2}$ denote the Lie product of
respective elements in $B$ and $A$.

\begin{lem} \label{lem:M'}
In this notation,
\begin{enumerate}
    \item the submodule $M_3= \left<Z_1 \cup Z_2 \right>$ of $M$ is
    isomorphic to $M_1 \oplus M_2$;
    \item the factor-module ${\raise0.7ex\hbox{$M$} \!\mathord{\left/
 {\vphantom {M_3}}\right.\kern-\nulldelimiterspace}
\!\lower0.7ex\hbox{$M_3$}}$ is isomorphic to $M_0$.
\end{enumerate}
\end{lem}

\proof To prove these statements we shall introduce the structure
of $R_{X \cup Y}$-module on the vector space $M'=M_0 \oplus M_1
\oplus M_2$ in such a way, that the vector space $M_1 \oplus M_2$
becomes a submodule of $M'$ and the induced action on $M_1 \oplus
M_2$ coincides with the initial one. Therefore, we only are to
introduce the action of $R_{X \cup Y}$ on the additive base of
$M_0$ in the correct way. Recall, that $M_0$ is isomorphic to the
ideal of $R_{X \cup Y}$ generated by the set $\left\{x_i y_j| i
\in I, j \in J \right\}$. Consequently, its additive base is  the
set of all monomials from elements of $X \cup Y$ which involve
letters from both of $X$ and $Y$.  Let $W$ be the set of all such
monomials written lexographically on $X$ and  on $Y$.

Let $w \in W$. Rewrite $w$ in the form \beq \label{eq:wW}
w=(x_iy_j)w^*,
 \eeq
 here $i$ and $j$ are the least indices of
letters form $X$ and $Y$ that are involved in $w$.

Let $V$ be the set of all monomials from $X \cup Y$. We define the
module action '$\cdot$' of an element $w \in W$ on $v \in V$ (and
denote it $w \cdot v$) using induction on the degree of $v$. We
set
\beq \label{eq:wx}
w\cdot x_{i_1}=\left\{
\begin{array}{ll}
    (x_iy_j)\cdot w^* x_{i_1}, & \hbox{if} \ i_1 \ge i; \\
    (x_{i_1}y_j)\cdot w^* x_i+(x_ix_{i_1})\cdot w^* y_j, & \hbox{if} \ i_1<i. \\
\end{array}
\right. \eeq Where $m=(x_ix_{i_1})\cdot w^* y_j$ is an element of
$M_1$ and $x_i x_{i_1}$ is the Lie product in the algebra $A$.

\beq \label{eq:wy} w\cdot y_{j_1}=\left\{
\begin{array}{ll}
    (x_iy_j)\cdot w^* y_{j_1}, & \hbox{if} \ j_1 \ge j; \\
    (x_iy_{j_1})\cdot w^* y_j+(y_{j_1}y_j)\cdot w^* x_i, & \hbox{if} \ j_1<j. \\
\end{array}
\right. \eeq Where $m=(y_{j_1}y_j)\cdot w^* x_i \in M_2$.

We leave the reader to check that Equations (\ref{eq:wx}) and
(\ref{eq:wy}) give rise to the required module structure on $M'$.

By the definition the module $M'$ is generated by the same
elements as the module $M$ is (the set coincide as words). We also
notice  that all the relations from Equation (\ref{eq:M}) which
are true in $M$ hold in $M'$. Therefore, there exists a
homomorphism $\phi: M \rightarrow M'$. We shall show that $\phi$
is an isomorphism. Consider the natural additive base
$\mathcal{B}$ for the module $M'$. The preimage of $\mathcal{B}$
under $\phi$ generates the module $M$. Moreover, all module
relations of $\mathcal{B}$ are satisfied in its preimage.
Consequently, $\phi(M_3)=M_1 \oplus M_2$ and the homomorphism
$\phi$ turns out to be an isomorphism and Lemma \ref{lem:M'} is
proven. \hfill $\blacksquare$

\medskip

By the module $M$ we construct a metabelian Lie algebra $C$. Let
$V$ be a vector space over the field $k$ with the base $X \cup Y$
and let $C=V \oplus M$ be the direct sum of vector spaces over
$k$. We next define the multiplication on the vector space $C$.

For a pair of elements $c_i \in C$, $c_i=(v_i,m_i); i=1,2; v_i \in
V, m_i \in M$ set
\begin{equation} \label{eq:mult}
(v_1,m_1)\circ (v_2,m_2)=(0,v_1\cdot v_2-m_2\cdot v_1+m_1\cdot
v_2).
\end{equation}
Here
$$
v_1v_2=\left\{%
\begin{array}{lll}
x_ix_j  & \hbox{ in } A, & \hbox{ if } v_1=x_i, v_2=x_j; \\
y_iy_j  & \hbox{ in } B, & \hbox{ if } v_1=y_i, v_2=y_j; \\
x_iy_j  & \hbox{ the element of } M, & \hbox{ if } v_1=x_i, v_2=y_j; \\
-x_iy_j & \hbox{ the element of } M, & \hbox{ if } v_1=y_i, v_2=x_j; \\
0,      &        & \hbox{ if } v_1=0  \hbox{ or } v_2=0.   \\
\end{array}%
\right.
$$
If $v_1$ and $v_2$ are linear combinations of elements from $X
\cup Y$, then the product $v_1v_2$ is defined using the
distributivity law and the equalities above.

Let $C$ be a vector space equipped with the operations introduced
above and $U(C)= \left\{(0,m)| m \in M \right\}$.

\begin{lem}
The set $U(C)$ is abelian ideal of the algebra $C$, $C^2=U(C)$ and
$C$ is a metabelian algebra.
\end{lem}
\proof Follows directly from the definitions of operations on $C$
\hfill $\blacksquare$

\medskip

\begin{lem}
The algebra $C$ is a Lie algebra, i.e. $C$ satisfies the
anti-commutativity identity and the  Jacoby identity. Therefore
$C$ is a metabelian Lie algebra.
\end{lem}
\proof 1. The anti-commutativity relation immediately follows from
(\ref{eq:mult})

2. Therefore, only the Jacoby identity is at issue. Consider a
3-tuple of elements from $C$. Let $x=(v_1,m_1)$, $y=(v_2,m_2)$,
$z=(v_3,m_3)$. By the definition
$$
x \circ y \circ z = (0,(v_1v_2)v_3+ (m_1v_2)v_3-(m_2v_1)v_3)
$$
and
$$
x \circ y \circ z + y \circ z \circ x + z \circ x \circ y =
(0,(v_1v_2)v_3+(v_2v_3)v_1+(v_3v_1)v_2).
$$
We now check that
\beq \label{eq:=0}
(v_1v_2)v_3+(v_2v_3)v_1+(v_3v_1)v_2=0.
\eeq

To establish equality (\ref{eq:=0}) we need to verify it on the
elements from $X \cup Y$ only. Suppose that $v_1, v_2$ and $v_3$
are the elements from either $X$ or $Y$. In which case
(\ref{eq:=0}) holds, since the Jacoby identity holds in both $A$
and $B$. Consider the case when $v_i; i=1,2,3$ lie in both $X$ and
$Y$. Due to the symmetry we may assume that $v_1=x_{i_1}$,
$v_2=x_{i_2}$, $v_3=y_j$. Expression (\ref{eq:=0}) takes the form
$$
(v_1v_2)v_3+(v_2v_3)v_1+(v_3v_1)v_2= (x_{i_1}x_{i_2})\cdot y_j+
(x_{i_2}y_j)\cdot x_{i_1} - (x_{i_1}y_j)\cdot x_{i_2}.
$$
There are three alternatives:
\begin{description}
    \item[If $x_{i_1}=x_{i_2}$] then Equation (\ref{eq:=0})
    obviously holds.
    \item[If $x_{i_1}>x_{i_2}$] then $(x_{i_1}y_j)\cdot x_{i_2}=
    (x_{i_1}x_{i_2})\cdot y_j+
    (x_{i_2}y_j)\cdot x_{i_1}$ and (\ref{eq:=0}) follows.
    \item[If $x_{i_2}>x_{i_1}$] analogous to the case above.
\end{description}
\hfill $\blacksquare$

\medskip

Consider the  algebra $A\ast B$. Clearly, it allows the generating
set $X \cup Y \cup Z_1 \cup Z_2$. Categorial Definition
\ref{def:A*B2} of the algebra $A\ast B$ yields to the existence of
the natural homomorphism $\psi: A\ast B \rightarrow C$.

\begin{thm} \label{thm:A*B1}
The homomorphism $\psi$ is an isomorphism. Moreover, $\psi ((A\ast
B)^2) \cong U(C) \cong M$. The commutant $A^2$ generates a
submodule of $(A\ast B)^2$ isomorphic to $M_1$ and $B^2$ generates
a submodule of $(A\ast B)^2$ isomorphic to $M_2$. The
factor-module $ {\raise0.7ex\hbox{$(A\ast B)^2$} \!\mathord{\left/
{\vphantom {M_1 \oplus M_2}}\right.\kern-\nulldelimiterspace}
\!\lower0.7ex\hbox{$(M_1 \oplus M_2)$}}$ is isomorphic to $M_0$.
\end{thm}
\proof Consider the natural additive basis of the algebra $C$ and
its natural preimage in $A\ast B$. Direct calculations show, that
this preimage in $A\ast B$ is a set of additive generators for it
(see \cite{Art}). Therefore, $\psi$ is an isomorphism. Now all
other statements of Theorem \ref{thm:A*B1} follow.
 \hfill $\blacksquare$

\medskip

\begin{thm} \label{thm:A*B2}
Let $A$ and $B$ be abelian Lie algebras over $k$ thus $A$ is a
vector space over $X$ and $B$ is a vector space over $Y$. Then
$R_{X \cup Y}$-module $(A\ast B)^2$ is isomorphic to $M_0$.
\end{thm}
\proof The conditions above imply, that $A = \overline{A}$ and
$B=\overline{B}$, therefore, Theorem \ref{thm:A*B2} is a direct
corollary of Theorem \ref{thm:A*B1} (in this cas $M_1=M_2=0$).
 \hfill $\blacksquare$

\medskip

\begin{prop} \label{prop:II-fit=comforprod}
The Fitting's radical of $A\ast B$, $A \ne 0, B \ne 0$ coincides
with the commutant $(A\ast B)^2$, {\rm
$$
\Fit(A\ast B)=(A\ast B)^2
$$
}
\end{prop}
\proof For every metabelian Lie algebra holds its Fittings radical
contains its commutant, thus $\Fit(A\ast B) \supseteq (A \ast
B)^2$. It is sufficient to show the inverse inclusion.  Assume the
converse $a \in \Fit(A\ast B)$ and $a \notin (A\ast B)^2$.
Consider $c=x_iy_j$, where $x_i \in X, y_j \in Y$ are
 chosen to be the least elements from $X$ and $Y$ respectively.
Let $a=v_1+v_2 +m$. Here $v_1$ is a linear combination of elements
from $X$, $v_2$ is a linear combination of elements from $Y$ and
$m \in (A\ast B)^2$. The assumption $a \notin (A\ast B)^2$ implies
that $c \circ \stackrel{n \ times}{\overbrace{a \circ \cdots \circ
a}} \neq 0$. Therefore, since  $a, ca \in \left< a \right>$ we
conclude that the ideal $\left< a \right>$ is not nilpotent and $a
\notin \Fit(A\ast B)$, what derives a contradiction.
 \hfill $\blacksquare$

\medskip

\section{Semidomains}

In this section we introduce the notions of a zero divisor, a
semidomain, a strict semidomain and establish the criterion which
points out when the metabelian product $A\ast B$ is a semidomain.

Let $A$ be Lie algebra over field $k$.

\begin{defn} \label{defn:0div}
A nonzero element \/ $x \in A$ is termed a zero divisor if and
only if there exists $y \in A, y \ne 0$ such that:
\begin{equation} \label{eq:0d1}
 \left< x \right> \circ \left< y \right> =0.
\end{equation}
\end{defn}

\medskip

Definition \ref{defn:0div} is symmetric on $x$ and $y$, thus we
say, that $x$ and $y$ is a pair of zero divisors.

In particular, if $A$ is a metabelian Lie algebra, then a pair $x,
y$ is a pair of zero divisors if and only if
\begin{equation} \label{eq:0d}
 xy=0; \ axy=0 \ \forall a\in A.
\end{equation}
Applying the Jacoby identity we see that (\ref{eq:0d}) is
symmetric on $x$ and $y$,
$$
ayx=xy a+axy.
$$

\begin{ex} \label{ex:1}
Let $A$ be a nilpotent Lie algebra of nilpotency class $n$. Every
element $x \in A$ is a zero divisor. Choose $y$ to be an element
from the center of $A$. Then $\verb"id"\left< y \right>$ is a
one-dimensional $k$-vector space with basis $\left\{ y\right\}$.
The pair $(x,y)$ is a pair of zero divisors.
\end{ex}

\medskip

\begin{ex} \label{ex:2}
Let $A$ be an arbitrary Lie algebra over $k$. Then every element
of the Fitting's radical {\rm $\Fit (A)$}  of the algebra $A$ is a
zero divisor. Consider {\rm $x \in \Fit (A)$}. Then the ideal
$\verb"id"\left< x \right>$ is a nilpotent ideal of $A$. Applying
the argument similar to the one of Example \ref{ex:1} the
statement follows.
\end{ex}

\medskip

\begin{ex} \label{ex:3}
Let $A$ be a metabelian Lie $k$-algebra. Suppose that $x, y \in A;
x,y \ne 0$ and $xy=0, x \in A^2$. Then $x,y$ is a pair of zero
divisors. Indeed, for every element $a$ of $A$ the product $ay \in
A^2$, therefore $a y x=0$ and Equation (\ref{eq:0d}) is satisfied.
\end{ex}

\medskip

Denote $D(A)$ the set of all zero-divisors of algebra $A$ with
$0$.

Let $D(A)$ be the set of all zero-divisors of the algebra $A$
together with $0$.

\begin{defn} \label{def:sd}
A Lie algebra $A$ is termed a semidomain if and only if {\rm
$$
D(A)=\Fit(A).
$$
}
\end{defn}

\medskip

On behalf of Example \ref{ex:2}, $D(A) \supseteq \Fit(A)$ holds
for every Lie algebra.

In paper \cite{AG1} we have introduced the notion of a
\emph{matrix metabelian Lie algebra}, which is constructed using a
free module over the rings of polynomials. One can check, that
every matrix metabelian Lie algebra is a semidomain.

\begin{ex}[Non-semidomain] \label{ex:4}
Consider a metabelian Lie algebra $A$ given by the following
generators and relations
$$
A= \left< a_1, a_2, a_3 |a_1a_2a_3=0 \right>_{\mathcal{A}^2}.
$$
Set $x=a_1a_2$, $y=a_3$. Notice that $x \in A^2$, $y \notin A^2$
and $xy=0$, which implies that $x$ and $y$ is a pair of zero
divisors. Consider non-zero product $a_3a_1a_3a_3 \cdots a_3 \ne
0$. Consequently, the ideal $\left< a_3 \right>$ is not nilpotent,
{\rm $y \notin \Fit(A)$} and thus $A$ is not a semidomain.
\end{ex}

\medskip

\begin{defn} \label{def:ssd}
A Lie algebra $A$ is termed a strict semidomain if and only if:
$$
D(A)=A^2.
$$
\end{defn}

\medskip

In the event that  $A$ is a metabelian Lie algebra we refine this
definition. For $\Fit(A) \supseteq A^2$, then $D(A) \supseteq
A^2$, which implies that $A$ is a strict semidomain if and only if
\begin{itemize}
    \item $A$ is a semidomain and
    \item $\Fit(A)= A^2$.
\end{itemize}

On behalf of Example \ref{ex:3}, if a metabelian Lie algebra is a
strict semidomain, then for all $x, y \in A; x,y \ne 0$
\begin{equation} \label{eq:ltf}
x \in A^2, y \notin A^2 \rightarrow x y \ne 0.
\end{equation}
In other words, $A^2$ regarded as module over the ring of
polynomials has no torsion under the action of linear polynomials.
In that case we  say that the  module is \emph{linear-torsion
free}.

\begin{thm}[Criterion of 'being strict semidomain'] \label{thm:cr}
Non-zero metabelian Lie algebra $A$ is a strict semidomain if and
only if:
\begin{itemize}
    \item $A$ is not abelian and
    \item $A^2$ is linear-torsion free.
\end{itemize}
 \end{thm}
\proof Let $A$ be a strict semidomain. Therefore, for all $x,y \in
A$ equation (\ref{eq:ltf}) holds. If $A$ is an abelian Lie
algebra, then $A^2=0$ and $D(A)=A$, thus, the case  when $A$ is an
ableian Lie algebra derives a contradiction.

Let $A$ be a nonabelian metabelian Lie algebra and (\ref{eq:ltf})
holds. We next derive, that $A$ is a strict semidomain. Consider a
pair of zero divisors $x, y$.

Suppose that  $x \in A^2$. Then (\ref{eq:ltf}) yields to $y \in
A^2$. Now suppose that $x \notin A^2$. For $A$ is not abelian
 and (\ref{eq:ltf}) holds there is to exist $z \in A, z \ne 0$
 such that $zx \ne 0$. On the
other hand, since $x,y$ is a pair of zero divisors we have that
$zxy=0$. Consequently, the pair $zx, y$ is a pair of zero
divisors. Applying the argument similar to the one above we
conclude that $y \in A^2$. Therefore $x \in A^2$ --- a
contradiction. \hfill $\blacksquare$

\medskip

\begin{cor} \label{cor:ssd}
The commutant $A^2$ of nonzero metabelian Lie algebra is
linear-torsion free if and only if $A$ is either abelian or a
strict semidomain.
\end{cor}

\medskip

Let $A\ast B$ be the metabelian product of metabelian Lie algebras
$A$ and $B$.  We have investigated the structure of this algebra
in Section \ref{sec:prod}, where we have shown that $\Fit(A \ast
B)=(A \ast B)^2$. This indicates that the properties of `being a
semidomain' and `being a strict semidomain' are either true or
false for $A\ast B$ simultaneously.

On behalf of Theorem \ref{thm:A*B1} the structure of $(A\ast B)^2$
is defined by three $R_{X \cup Y}$-modules $M_0$, $M_1$ and $M_2$.
The module $M_0$ is torsion free, while in general $M_1$ and $M_2$
are not (for instance, they contain $A^2$ and $B^2$ as
submodules). Although the following lemma holds:

\begin{lem} \label{lem:as}
Assume that the commutant $A^2$ of the algebra $A$ is
linear-torsion free. Then the module $M_1$ is linear-torsion free.
\end{lem}
\proof Consider $u \in M_1$ and a linear polynomial $f \in R_{X
\cup Y}$:
$$
f= \sum \alpha_i x_i + \sum \beta_j y_j.
$$
Take the  summand with the greatest singleton in $u$, such that it
involves the letters of $Y$ only:
$$
u=u_1+a \otimes y_{j_1}^{g_1} \cdots y_{j_t}^{g_t}, \ 0\ne a \in
A^2
$$
Suppose that $\sum \beta_j y_j \ne 0$ and $y_{j_n}$ is the letter
with the greatest index and non-zero coefficient among the letters
that occur in the linear combination $\sum \beta_j y_j$. Then in
the product $u \cdot f$ the summand $\beta_j a \otimes
y_{j_1}^{g_1} \cdots y_{j_t}^{g_t} \cdot y_{j_n}$ is not equal to
zero and does not cancel.

Now suppose that $\sum \beta_j y_j = 0$. This implies that $\sum
\alpha_i x_i \ne 0$. According to the restrictions of Lemma
\ref{lem:as} $A^2$ is linear torsion free and we obtain:
$$
a \cdot \sum \alpha_i x_i \ne 0.
$$
Therefore, in the product $u\cdot f$ the summand
$$
(a \cdot \sum \alpha_i x_i) \otimes y_{j_1}^{g_1} \cdots
y_{j_t}^{g_t}
$$
does not cancel.

Thus, $u \cdot f \ne 0$ --- q. e. d.
 \hfill $\blacksquare$

\medskip

Finally we formulate a criterion for $A\ast B$ to be a semidomain.

\begin{thm} \label{thm:stA*B}
Let $A$ and $B$ be two non-zero metabelian Lie algebras over $k$.
Then $A\ast B$ is a semidomain if and only if each of $A$ and $B$
is either abelian or a strict semidomain.
\end{thm}
\proof As mentioned above, the properties of 'being a semidomain'
and 'being a strict semidomain' are equivalent for $A\ast B$.
Since $A$ and $B$ are non-zero, the algebra $A\ast B$ is
non-abelian. Theorem \ref{thm:cr} states that $A\ast B$ is a
semidomain if and only if its commutant $(A\ast B)^2$ is
linear-torsion free.

Assume that $A$ is neither abelian nor strict semodomain. Then the
commutant $A^2$ has linear torsion, i.e. there are non-zero
elements $x$ and $y$ of $A$ such that $x \in A^2$, $y \notin A^2$
and $xy=0$. Therefore, $x \in (A\ast B)^2$, $y \notin (A \ast
B)^2$ what essentially implies that $(A\ast B)^2$ has linear
torsion. Due to Theorem \ref{thm:cr} $A \ast B$ is not a
semidomain.

Now assume, that each of $A$ and $B$ is either abelian or a strict
semidomain. Thus, the commutants $A^2$ and $B^2$ are linear
torsion free. We are to show that $(A \ast B)^2$ is linear torsion
free.

Consider $c \in (A \ast B)^2, c \ne 0$ and a linear polynomial $f
\in R_{X \cup Y}$:
$$
f= \sum \alpha_i x_i + \sum \beta_j y_j.
$$
It is sufficient to show, that $c \cdot f \ne 0$.

Since $M_0$ is a torsion free module $c \notin M_1 \oplus M_2$ we
obtain that $c \cdot f \ne 0$. Consequently, let $c \in M_1 \oplus
M_2$ and write it as follows:
$$
c=u+v, \ u\in M_1,\ v \in M_2.
$$
Lemma \ref{lem:as} states, that either $u\cdot f \ne 0$ or $v
\cdot f \ne 0$. For $M_1 \cap M_2 = \emptyset$ we have $u\cdot f +
v\cdot f \ne 0$, i.e. $c \cdot f \ne 0$.
 \hfill $\blacksquare$

\textit{Evelina Yu. Daniyarova, e-mail:
\underline{evelina.omsk@list.ru}
\newline
644099, Russia, Omsk, Pevtsova st. 13, Omsk Branch of Institute of
Mathematics (Siberian Branch of Russian Academy of Science).}
\newline
Telephone: +7 381-2-23-67-39

\medskip

\textit{Ilia V. Kazatchkov, e-mail:
\underline{kazatchkov@iitam.omsk.net.ru}
\newline
644099, Russia, Omsk, Pevtsova st. 13, Omsk Branch of Institute of
Mathematics (Siberian Branch of Russian Academy of Science).}
\newline
Telephone: +7 381-2-23-67-39

\medskip

\textit{Vladimir N. Remeslennikov, e-mail:
\underline{remesl@iitam.omsk.net.ru}
\newline
644099, Russia, Omsk, Pevtsova st. 13, Omsk Branch of Institute of
Mathematics (Siberian Branch of Russian Academy of Science).}
\newline
Telephone: +7 381-2-23-67-39
\end{document}